\documentclass[conference]{IEEEtran}
\IEEEoverridecommandlockouts
\usepackage{cite}
\usepackage{xparse, xspace}
\usepackage{amsmath, amsfonts, mathtools, stmaryrd, amsbsy, bm, breqn, amsthm}
\usepackage[mathscr]{eucal}
\usepackage{algorithm}
\usepackage{algorithmicx,algpseudocode}
\usepackage{graphicx}
\usepackage{textcomp}
\usepackage{xcolor}
\usepackage{fullpage, graphicx, hyperref}
\usepackage{tikz}
\usetikzlibrary{external}
\usepackage{float}
\usepackage{cleveref}
\usepackage{multirow}

\usepackage{graphicx}
\usepackage{xcolor}
\newcommand{\stirling}[2]{\genfrac{\{}{\}}{0pt}{}{#1}{#2}}

\newcommand{\ttsv}[1][]{\textsc{TTSV#1}}

\usepackage{braket}
\usepackage{siunitx}
\NewDocumentCommand{\Tn}{O{} m}{\boldsymbol{#1\mathscr{\MakeUppercase{#2}}}}

\NewDocumentCommand \Vc { O{} m } {{\bm{#1\mathbf{\MakeLowercase{#2}}}}}
\NewDocumentCommand \Mx { O{} m } {{\bm{#1\mathbf{\MakeUppercase{#2}}}}} 

\Crefname{equation}{Eq.}{Eqs.}
\Crefname{figure}{Fig.}{Figs.}
\Crefname{tabular}{Tab.}{Tabs.}
\Crefname{section}{Sec.}{Secs.}
\crefname{thm}{theorem}{theorems}
\Crefname{thm}{Theorem}{Theorems}
\crefname{prop}{property}{properties}
\Crefname{prop}{Property}{Properties}
\usepackage{ourmacros}
\def\BibTeX{{\rm B\kern-.05em{\sc i\kern-.025em b}\kern-.08em
    T\kern-.1667em\lower.7ex\hbox{E}\kern-.125emX}}

\makeatletter
\newcommand{\linebreakand}{%
  \end{@IEEEauthorhalign}
  \hfill\mbox{}\par
  \mbox{}\hfill\begin{@IEEEauthorhalign}
}
\makeatother
\begin{document}

\title{Fast Parallel Tensor Times Same Vector for Hypergraphs
\thanks{Authors IA, SGA, JL, and SJY gratefully acknowledge the funding support from the Applied Mathematics Program within the U.S.\ Department of Energy’s Office of Advanced Scientific Computing Research as part of Scalable Hypergraph Analytics via Random Walk Kernels (SHARWK). Pacific Northwest National Laboratory is operated by Battelle for the DOE under Contract DE-AC05-76RL0 1830. PNNL Information Release PNNL-SA-187496}
}

\author{
\IEEEauthorblockN{Shruti Shivakumar}
\IEEEauthorblockA{School of Computational \\
Science and Engineering\\
Georgia Institute of Technology\\
Atlanta, GA \\
sshivakumar9@gatech.edu}\\   
\IEEEauthorblockN{Jiajia Li}
\IEEEauthorblockA{School of Computer Science\\
North Carolina State University\\
Raleigh, NC\\
jiajia.li@ncsu.edu}
\and
\IEEEauthorblockN{Ilya Amburg}
\IEEEauthorblockA{Pacific Northwest\\
National Laboratory\\
Richland, WA\\
ilya.amburg@pnnl.gov}\\[0.4cm]  
\IEEEauthorblockN{Stephen J.\ Young}
\IEEEauthorblockA{Pacific Northwest\\
National Laboratory\\
Richland, WA\\
stephen.young@pnnl.gov}
\and
\IEEEauthorblockN{Sinan G.\ Aksoy}
\IEEEauthorblockA{Pacific Northwest \\
National Laboratory\\
Seattle, WA\\
sinan.aksoy@pnnl.gov}\\ [0.4cm]                
\IEEEauthorblockN{Srinivas Aluru}
\IEEEauthorblockA{School of Computational \\
Science and Engineering\\
Georgia Institute of Technology\\
Atlanta, GA \\
aluru@cc.gatech.edu}
}


\maketitle

\begin{abstract}
Hypergraphs are a popular paradigm to represent complex real-world networks exhibiting multi-way relationships of varying sizes. Mining centrality in hypergraphs via symmetric adjacency tensors has only recently become computationally feasible for large and complex datasets. To enable scalable computation of these and related hypergraph analytics, here we focus on the \underline{S}parse \underline{S}ymmetric \underline{T}ensor \underline{T}imes \underline{S}ame \underline{V}ector (\sttsvc{}) operation. We introduce the \cmpdformatfull{} (\cmpdformat{}) format, an extension of the compact \format{} format for hypergraphs of varying hyperedge sizes and present a shared-memory parallel algorithm to compute \sttsvc{}. We experimentally show \sttsvc{} computation using the \cmpdformat{} format achieves better performance than the naive baseline, and is subsequently more performant for hypergraph $H$-eigenvector centrality.
\end{abstract}



\section{Introduction} \label{sec:introduction}

Hypergraphs are generalizations of graphs that represent multi-entity relationships in a broad range of domains, such as cybersecurity \cite{Yamaguchi2015CyberAlgorithms,Joslyn2020HypergraphRelationships}, biological systems \cite{Feng2021HypergraphResponse}, social networks \cite{Zhu2019SocialNetworks}, and telecommunications\cite{Ganesan2023StructuredSystems}. 
While graph edges connect exactly two nodes, hyperedges may connect any number of nodes. Most real-world hypergraphs are non-uniform, meaning they have differently sized hyperedges, which poses challenges for their compact representation and analysis. 

For uniform hypergraphs, symmetric tensors are popularly used to represent the higher-order adjacency information \cite{Benson2015TensorStructures,Ke2020CommunityIteration}. Several strategies have been explored to extend the tensor representation approach to non-uniform hypergraphs, including adding dummy nodes \cite{Ouvrard2018OnTensor,Zhen2021CommunityEmbedding}, considering the set of symmetric adjacency tensors, with each tensor arising from the component uniform hypergraph of a particular hyperedge size \cite{Dumitriu2021PartialModel,Dumitriu2023ExactModel}, and combinatorially inflating the lower--cardinality hyperedges until all hyperedges are equisized \cite{Banerjee2017SpectraHypergraphs}. 
Following Aksoy, Amburg, and Young \cite{Aksoy2023ScalableHypergraphs}, we call these inflated hyperedges \textit{blowups} and focus on the adjacency tensor as defined by Banerjee et al. \cite{Banerjee2017SpectraHypergraphs}, which we call the \textit{\blowup{}} tensor.

\sttsvc{} is a key operation on symmetric tensors, and is the computational bottleneck in fundamental algorithms such as the shifted-power method for computing tensor eigenpairs and symmetric CP-decomposition \cite{Kolda2015NumericalDecomposition,Kolda2011ShiftedEigenpairs,Kolda2015SymmetricTrivial}. These algorithms, in turn, are utilized to perform a variety of hypergraph analyses. 
For instance, eigenpairs of the adjacency tensor of uniform hypergraphs are used to define $H$-eigenvector centrality (HEC) \cite{Benson2019ThreeCentralities}, a nonlinear hypergraph centrality measure which was further extended to non-uniform hypergraphs using the \blowup{} tensor representation \cite{Aksoy2023ScalableHypergraphs}.
Thus, developing performant algorithms for \sttsvc{} on the \blowup{} tensor enables efficient computation of hypergraph centrality.   

However, working with the \blowup{} tensor requires we overcome several computational challenges. 
First, since enumerating all its nonzeros is prohibitively costly, following Aksoy, Amburg, and Young, we will adapt the ``generating function approach" \cite{Aksoy2023ScalableHypergraphs}, to perform the computation {\it indirectly}.
Second, we introduce a new, compressed format for tensors that is tailored to reduce the memory footprint of storing nonuniform hypergraphs, called \cmpdformatfull{} (\cmpdformat{}). This extends past work on the \format{} format for uniform hypergraphs\cite{Shivakumar2021EfficientTensors,Shivakumar2023SparseDecomposition} and, as explained further in Section \ref{sec:s3ttvc}, achieves \sttsvc{} performance gains via memoization of intermediate results.

Our main contributions are summarized as follows:
\begin{itemize}
    \item We introduce the \cmpdformatfull{} (\cmpdformat{}), an extension of the \format{} format for non-uniform hypergraphs, and demonstrate up to $26.4\times$ compression compared to coordinate storage format for real-world hypergraphs.
    \item We implement an efficient multi-core parallel \sttsvc{} algorithm, called \algo{} which adapts the generating function approach to the \cmpdformat{} format, and identifies opportunities for memoization of intermediate results. 
    \item We present two baseline approaches which use the \cmpdformat{} without memoization and adopt two state-of-the-art approaches to highlight the performance of \algo{}. We realize up to $53.98\times$ speedup compared to \naive{}, and up to $12.45\times$ speedup compared to \fft{}. 
    \item We apply our algorithm to the calculation of $H$-eigenvector centrality for hypergraphs, obtaining speedups of many orders of magnitude over state-of-the-art approaches.
\end{itemize}
\section{Preliminaries} \label{sec:preliminaries}
Following Kolda and Bader~\cite{Kolda2009TensorApplications}, we denote vectors using bold lowercase letter (e.g., \V{a}, \V{b}),
and tensors using bold calligraphic letters (e.g., $\T{X}$).  For a tensor $\T{X}$ and a vector $\Vc{b}$ we will also denote by $\T{X} \times_j \Vc{b}$ the product of $\T{X}$ and $\Vc{b}$ along the $j^{\textrm{th}}$-mode of $\T{X}$ resulting in a order $(N-1)$ tensor. 
Using this notation, we can define 
Tensor-Time-Same-Vector in all modes but 1 (TTSV1) 
 \begin{equation}\label{eq:s3ttvc}    \V{s} = \T{X}\V{b}^{N-1} = \T{X} 
    = \T{X}  \times_2 \V{b} \times_3 \V{b} \ldots \times_N \V{b}
\end{equation}
which plays an important role in calculating generalized eigenvalues and eigenvectors associated with $\T{X}.$  Alternatively, by expanding along indicies this may be rewritten as
\begin{equation}\label{eq:ttsv1}
s_{i_1} = \left[\T{X}\Vc{b}^{N-1}\right]_{i_1} = \sum_{i_2=1}^n\cdots\sum_{i_N=1}^n\T{X}_{i_1,\dots,i_N}\prod_{k=2}^N\Vc{b}_{i_k}.
\end{equation}

Formally, a hypergraph is a pair $H=(V, E)$, where $V$ is the set of vertices, and $H$ is the set of hyperedges on those vertices; that is, $E$ is some subset of the power set of $V$. We say $H$ is uniform if all hyperedges have the same size; otherwise, it is nonuniform. The rank of $H$ is the size of the largest hyperedge. 

An order-$N$ symmetric tensor $\T{X}$ has $N$ modes or dimensions, with the special property that the values, $\T{X}_{(i_1, i_2, \ldots i_N)}$ remains unchanged under any permutation of its indices. Symmetric tensors arise naturally in many context including the representation of hypergraphs.  For example, if $H$ is an $N$-uniform hypergraph on $n$ vertices, there is natural symmetric representation as an order-$N$ tensor $\T{X} \in \Rn{n \times n \times \cdots \times n}$.  That is, for every hyperedge $e = \{v_{i_1}, v_{i_2}, \ldots, v_{i_N} \} \in E$ with weight $w(e)$, $\T{X}_{\sigma(\V{i})} = \frac{w(e)}{N!}$, where $\sigma(\V{i})$ denotes any of the $N!$ permutations of the index tuple $\V{i} = (i_1, \ldots, i_N)$.  In order to extend this representation to non-uniform hypergraphs, we follow the approach of Banerjee et al.~\cite{Banerjee2017SpectraHypergraphs} and define for a rank-$N$ edge-weighted hypergraph $H = (V,E,w)$ the order-$N$ \blowup{} tensor, $\blowupT$, associated with $H$.  To this end, for each edge $e \in H$ define the set of ordered blowups of $e$ as  
\[  \beta(e) = \cbr*{i_1, i_2, \ldots i_{N} : \text{for each } v \in e, \exists j \ni i_j = v}.\footnote{Note that $\lvert \beta(e)\rvert$ depends only on $N$ and the size of $e$ and is given by $\lvert e\rvert!\stirling{N}{|e|},$ where $\stirling{N}{|e|},$ is the Stirling number of the second kind.  Thus $\beta(e)$ can be trivially precomputed and stored in a lookup table to accelerate future computations.} \]
Then for each $\Vc{i} \in \beta(e)$, $\blowupT_{\Vc{i}}$ has value $\frac{w(e)}{|\beta(e)|}$. As we be working primarily with the blowup tensor, it will be convenient to define $\mathcal{E}(\blowupT)$ as the collection of edges which generated the blowup tensor $\blowupT$.  In this paper, we take $w(e) = \lvert e \rvert$ to ensure that $\blowupT\V{1}^{N-1} = \V{d},$ the vector of node degrees. It is worth noting that for uniform hypergraphs, $\blowupT$ is precisely the uniform adjacency tensor of the hypergraph discussed above.

The $H$-eigenvector centrality vector of $H$ is a positive vector $\V{x}$ satisfying $\blowupT\Vc{x}^{N-1}=\lambda\Vc{x}^{[N-1]},$ where $\lambda$ is the largest $H$-eigenvalue of $\blowupT$ and the vector operation $\V{x}^{[N-1]}$ represents componentwise $N-1$ power of $\V{x}$. By the Perron-Frobenius theorem for the hypergraph adjacency tensor~\cite{Aksoy2023ScalableHypergraphs}, if $H$ is connected, then $\V{x}$ is guaranteed to exist, and is unique up to scaling. Intuitively, here a node's importance (to the power of $N-1$, which guarantees dimensionality preservation) is proportional to a product of centralities over all blowups of hyperedges that contain it. A popular approach is to compute the eigenpair $(\lambda, \V{X})$ using the NQZ algorithm~\cite{Ng2010FindingTensor} (Algorithm~\ref{alg:nqz}, where $\oslash$ denotes componentwise division), an iterative power-like method that utilizes TTSV1 as its workhorse subroutine, which we employ here.
\begin{algorithm}
\footnotesize
    \begin{algorithmic}[1]
        \State \textbf{Input}: $n$-vertex, rank $N$ hypergraph $H$, tolerance $\tau$
        \State \textbf{Output}: $H$-eigenvector centrality, $\Vc{x}$
        \State $\V{y} = \frac{1}{n}\cdot \V{1}$
        \State $\V{z} = \ttsv[1](H, \V{y})$
        \Repeat
        \State $\Vc{x}=\Vc{z}^{\left[\frac{1}{N-1}\right]} / \norm{\Vc{z}^{\left[\frac{1}{N-1}\right]}}_1$
        \State $\Vc{z}=\ttsv[1](H,\Vc{x})$
        \State $\lambda_{\min} = \min{(\Vc{z} \oslash \Vc{x}^{\left[N-1\right]})}$
        \State $\lambda_{\max} = \max{(\Vc{z} \oslash \Vc{x}^{\left[N-
        1\right]})}$
        \Until{$(\lambda_{\max} -\lambda_{\min})/\lambda_{\min} < \tau$}
        \State \bf{return} $\V{x}$
    \end{algorithmic}
    \caption{NQZ algorithm for computing HEC}
    \label{alg:nqz}
\end{algorithm}

Motivated by questions in hypergraph node ranking, we investigate the TTSV1 operation for the blowup tensor of a non-uniform hypergraph.  To distinguish from the more general case, and emphasize the applicability to sparse symmetric tensors, we will refer to this problem as he \underline{S}parse \underline{S}ymmetric \underline{T}ensor \underline{T}imes \underline{S}ame \underline{V}ector (\sttsvc{}) operation on the \blowup{} tensor. 
In many ways, the current work can be thought of as synthesis of the implicit \sttsvc{} algorithm on the \blowup{} tensor proposed by Aksoy, Amburg and Young \cite{Aksoy2023ScalableHypergraphs}, with the \format{} format for storing sparse symmetric adjacency tensors of uniform hypergraphs \cite{Shivakumar2021EfficientTensors}.

To that end, we summarize some of the key features of these two approaches in the next two subsections.

\subsection{Generating Functions for \sttsvc{}}
Aksoy, Amburg and Young proposed the implicit AAY algorithm \cite{Aksoy2023ScalableHypergraphs} (\Cref{alg:without-css}) to evaluate TTSV1 for the \blowup{} tensor that relies on generating functions.
The fundamental observation which drives their algorithm is that, in the \blowup{} tensor, all entries corresponding to a single edge have the same coefficient.  Thus, by using generating functions to aggregate over the contributions of all elements of $\beta(e)$, the computational requirements can be significantly reduced.  More concretely, they observed that for any edge $e \in E$ and vertex $v \in e$, the contribution of $e$ to $[ \blowupT \V{b}^{N-1}]_v$ can be captured as a rescaling of the last entry in 
\[ E_N(b_v) \conv \left(\bigconv_{u \in e\backslash v}\overline{E}_r(b_u)\right), \]
where 
\begin{align*} 
E_N(c) &= \left[1, c, \frac{c^2}{2!}, \ldots, \frac{c^{N-1}}{(N-1)!} \right]  \quad \textrm{and} \\
\overline{E}_N(c) &= \left[0,c, \frac{c^2}{2!}, \ldots, \frac{c^{N-1}}{(N-1)!}\right]
\end{align*}
and $(a \conv b)$ is a vector of length $N+1$ representing the convolution operation with $(a \conv b)[k] = \sum_{i=0}^{k} a_i b_{k - i}$.  Alternatively, their approach 
can be viewed as extracting a specific coefficient of $t^{N-1}$ from a particular exponential generating function \cite{Wilf2006Generatingfunctionology}.  This approach yields \Cref{alg:without-css}. 
\begin{algorithm}
\footnotesize
    \begin{algorithmic}[1]
        \State \textbf{Input}: rank $N$ weighted hypergraph $(V,E,w)$, vector $\V{b}$
        \State \textbf{Output}: \sttsvc{} output, $\V{s} = \blowupT \V{b}^{N - 1}$
        \For{$v\in V$}

                    \State $c \gets 0$
                    \For{$e\in E(v)$}
                        \State \parbox[t]{2in}{$c\pluseq \displaystyle \frac{w(e)}{|\beta(e)|} (N-1)!E_N(\Vc{b}_v) \ast ( \bigast_{ u \in e\backslash v}\overline{E}_N(\Vc{b}_u))[N-1]$}
                    \EndFor
        \State $\V{s}_v \gets c$
        \EndFor
        \State \bf{return} $\V{s}$
    \end{algorithmic}
    \caption{AAY algorithm for implicit TTSV1 using Banerjee adjacency tensor. }
    \label{alg:without-css}
\end{algorithm}

While the aggregation over vertex-edge pairs given by the AAY approach results in significant computational speedups, the lack of structure imposed on the computation  results in frequent repetition of the convolution calculations. For example, in the AAY approach the convolution $\overline{E}(\Vc{b}_v) \conv \overline{E}(\Vc{b}_u)$ is computed $|e|-2$ times for every edge containing both $v$ and $u$.

\subsection{Compressed Sparse Symmetric Format}

The \format{} structure \cite{Shivakumar2021EfficientTensors,Shivakumar2023SparseDecomposition} is a compact storage format that enables efficient \sttsmc{} computation for sparse symmetric adjacency tensors $\T{X}$ arising from uniform hypergraphs.  In order to take advantage of the symmetry of $\T{X}$,  \format{} stores all information based on the collection of sorted edges of the associated hypergraph, $\mathcal{E}(\T{X})$\footnote{Vertex-sorted hyperedges are also referred to as \emph{index-ordered unique (IOU) non-zeros} \cite{Shivakumar2021EfficientTensors,Shivakumar2023SparseDecomposition}, and denoted by $unz(\T{X})$, where a non-zero entry of a order $N$ tensor $\T{X}_{i_1,\ldots, i_N}$ is IOU if $i_1 < i_2 < \cdots < i_N$.}

If $\T{X}$ has order $N$, then the \format{} is a forest with $N-1$ levels where every length $k$ subsequence of an element of $\mathcal{E}(\T{X})$ corresponds to a unique root to level $k$ path in the \format{}.  Further, the leaves at level $N-1$ are equipped with the "dropped" index in the subsequence and the value of $\T{X}$ for the corresponding element of $\mathcal{E}(\T{X})$.  A key advantage of using the \format{} format is its computation-aware nature --- by storing all ordered subsequences of $\mathcal{E}(\T{X})$, intermediate results in the \sttsmc{} computation can be easily memoized with minimal additional index information.
In Shivakumar et al. \cite{Shivakumar2021EfficientTensors, Shivakumar2023SparseDecomposition} the \sttsmc{}-\format{} algorithm is used to find the tensor decomposition of an adjacency tensor of a uniform hypergraph. The convergence of this method requires that the original hypergraph be connected. For a non-uniform hypergraph, it is likely that there exists an edge size such that the collection of hyperedges of that size is not connected. Thus, it is theoretically necessary to work with a single tensor representation of the hypergraph, such as the blowup tensor, in order to preserve the necessary convergence properties. This presents two primary challenges in applying \sttsmc{}-\format{} that the current work addresses: the blowup tensor can have super-exponentially many non-zeros corresponding to a single edge and any data structure must explicitly account for the repetitions of the vertices induced by the blow-up. Naively, extending the index-ordered non-zeros approach of \sttsmc{}-\format{} to incorporate repeated vertices will result in a significant increase in the memory footprint of the \format{} structure to account for the repeated vertices, as well as the computational cost of computing \sttsvc{} itself. 
On the other hand, adopting the implicit construction approach and storing the adjacency tensor of each constituent uniform hypergraph using the \format{} format is a suboptimal approach in terms of memory requirement compared to \cmpdformat{} (described in the next section) and results in greater computation cost as we lose out on memoizing intermediate $\bar{E}_N$ across IOU non-zeros i.e. hyperedges.
Moreover, directly applying the \sttsmc{}-\format{} algorithm to such a storage construction
will not lead to correct results for the tensor-times-same-vector operation. 

\section{Compound CSS Structure}

In this paper, we present an extension of \format{}, called \cmpdformatfull{} (\cmpdformat{}) which facilitates fast \sttsvc{} computation on the \blowup{} tensor representing non-uniform hypergraphs.  
One natural way to extend \format{} to non-uniform hypergraphs would be to build an instance of \format{} for 
each constituent uniform hypergraph and work independently on the corresponding adjacency tensors.
However, as the associated tenors would all be of different orders, additional work would be needed to lift the methods of Shivakumar et al.~\cite{Shivakumar2021EfficientTensors} to the non-uniform case.  Furthermore, similar to the AAY approach, decomposing the non-uniform hypergraph into multiple uniform hypergraphs leads to significant extra computation and storage. For instance, if $u$ and $v$ are in multiple edges of different sizes then memoized work for the subsequence $u, v$ occurs in the \format{} for each of these edge sizes. 
 
 Instead, we introduce \cmpdformat{} which extends the \format{} to non-uniform hypergraphs by building a forest of $N-1$ levels containing all proper subsequences $\mathcal{E}(\T{B})$, that is, the ordered proper subsets of the edges.  In particular, if $f$ is a size $\ell$ proper subset of an edge $e$, $f$ is represented by a unique root to level $\ell$ path in the forest, and further, this path is given by an in-order listing of the elements of $f$.  In contrast to \format{}, in \cmpdformat{} the edges of $\T{B}$ are ``owned" by vertices in the forest at all levels and a vertex at a given level may own multiple edges.  We note that the edges owned by a given vertex can be thought of as special leaves of the data structure (at level corresponding the the edge size) which store the ``dropped" vertex and the value of the tensor at all blowups of the edge.  These special leaves of the \cmpdformat{} can be easily enumerated as the ordered pairs $\mathcal{L} = \{ (e,v) : e \in E, v \in e\}$.  We will denote by $\mathcal{L}_k \subseteq \mathcal{L}$ those special leaves corresponding to an edge of size $k.$  We will also denote by $\mathcal{S}(v)$ the set of special leaves ``owned" by a vertex in the \cmpdformat{} structure and note that $\mathcal{L}_{\ell} = \cup \mathcal{S}(v)$ where the union is taken over all vertices at level $\ell-1$ in the \cmpdformat{} structure.
 
\begin{figure}
\begin{minipage}{0.5\linewidth}
\centering
\begin{minipage}[b]{.1\linewidth}
\centering
\raisebox{-.45\height}{\includegraphics[width=40pt]{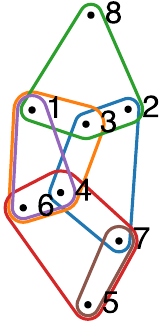}}
\end{minipage}
\hspace{35pt} 
\begin{minipage}[b]{0.35\linewidth} 
\centering
    \scalebox{1}{
    \begin{tabular}{c | c c c c c c c}
        $i_1$ & 2 & 1 & 1 & 4 & 1 & 5\\
        $i_2$ & 3 & 3 & 2 & 5 & 4 & 7\\
        $i_3$ & 4 & 4 & 3 & 6 & 6\\
        $i_4$ & 7 & 6 & 8 & 7\\
        \hline
        vals & $\color{blue}v_1$ & $\color{orange}v_2$ & $\color{darkgreen}v_3$ & $\color{red}v_4$ & $\color{violet}v_5$ & $\color{brown}v_6$
    \end{tabular}
    }
\end{minipage}
\end{minipage}
\caption{Non-uniform hypergraph generating three sparse symmetric tensors having $6$ IOU nonzeros in total.}\label{tab:example}
\end{figure}

In the example shown in \Cref{fig:ccss}, the \cmpdformat{} is constructed from a 8-node weighted non-uniform hypergraph with edges (shown in index-ordered format) in \Cref{tab:example}. We can easily see the reduction in space (as compared to the multiple \format{} approach) in this example -- the sequence $(1,4)$ is shared between the sequences $(1,4)$ and $(1,4,6)$, corresponding to the edges $\{1,4,6\}$ and $\{1,3,4,6\}$, respectively.

\begin{figure*}
\centering
	\begin{tikzpicture}[level distance=10mm]
		\tikzstyle{regular}=[circle, draw, inner sep=1pt, font=\footnotesize]
		\tikzstyle{level 1}=[sibling distance=25mm]
		\tikzstyle{level 2}=[sibling distance=15mm]
		\tikzstyle{level 3}=[sibling distance=10mm]
		\tikzstyle{level 4}=[sibling distance=7mm]
		\tikzstyle{phantom}=[]
		\tikzstyle{val}=[rectangle, align=center, inner sep=1pt]
		\node [phantom] {}
		[every node/.style={regular}]
		child {node {1} edge from parent[draw=none]
			child [xshift=14mm]{node {2}
				child[xshift=2mm]{node {3}
					child {node [val] {8\\$\color{darkgreen}v_3$} edge from parent[dashed]}
				}
				child [xshift=-3mm]{node {8}
					child {node [val] {3\\$\color{darkgreen}v_3$} edge from parent[dashed]}
				}
			}
			child [xshift=11mm]{node {3}
				child[xshift=5mm]{node {4}
					child {node [val] {6\\$\color{orange}v_2$} edge from parent[dashed]}
				}
				child[xshift=0mm]{node {6}
					child {node [val] {4\\$\color{orange}v_2$} edge from parent[dashed]}
				}
				child[xshift=-5mm]{node {8}
					child {node [val] {2\\$\color{darkgreen}v_3$} edge from parent[dashed]}
				}
			}
			child [xshift=8mm]{node {4}
				child[xshift=3mm]{node {6}
					child {node [val] {3\\$\color{orange}v_2$} edge from parent[dashed]}
				}
                child [xshift=-2mm]{node [val] {6\\$\color{violet}v_5$} edge from parent[dashed]
                }
			}
			child [xshift=2mm]{node {6}
                child{node [val] {4\\$\color{violet}v_5$} edge from parent[dashed]
                }
            }
			child [xshift=-5mm]{node {8}}
		}
		child [xshift=17mm]{node {2} edge from parent[draw=none]
			child [xshift=12mm]{node {3}
				child [xshift=5mm]{node {4}
					child {node [val] {7\\$\color{blue}v_1$} edge from parent[dashed]}
				}
				child [xshift=0mm] {node {7}
					child {node [val] {4\\$\color{blue}v_1$} edge from parent[dashed]}
				}
				child[xshift=-5mm]{node {8}
					child {node [val] {1\\$\color{darkgreen}v_3$} edge from parent[dashed]}
				}
			}
			child[xshift=8mm]{node {4}
				child{node {7}
					child {node [val] {3\\$\color{blue}v_1$} edge from parent[dashed]}
				}
			}
			child [xshift=-1mm]{node {7}}
			child [xshift=-10mm]{node {8}}
		}
		child [xshift=22mm]{node {3} edge from parent[draw=none]
			child [xshift=11mm]{node {4}
				child[xshift=2mm]{node {6}
					child {node [val] {1\\$\color{orange}v_2$} edge from parent[dashed]}
				}
				child[xshift=-3mm]{node {7}
					child {node [val] {2\\$\color{blue}v_1$} edge from parent[dashed]}
				}
			}
			child [xshift=2mm]{node {6}}
			child [xshift=-8mm]{node {7}}
			child [xshift=-18mm]{node {8}}
		}
		child[xshift=12mm]{node {4} edge from parent[draw=none]
			child [xshift=10mm]{node {5}
				child [xshift=3mm]{node {6}
					child {node [val] {7\\$\color{red}v_4$} edge from parent[dashed]}
				}
				child[xshift=-2mm]{node {7}
					child {node [val] {6\\$\color{red}v_4$} edge from parent[dashed]}
				}
			}
			child [xshift=5mm]{node {6}
				child[xshift=3mm]{node {7}
					child {node [val] {5\\$\color{red}v_4$} edge from parent[dashed]}
				}
                child{node [val] {1\\$\color{violet}v_5$} edge from parent[dashed]
                }
			}
			child [xshift=-3mm]{node {7}}
		}
		child [xshift=10mm] {node {5} edge from parent[draw=none]
			child [xshift=10mm]{node {6}
				child{node {7}
					child {node [val] {4\\$\color{red}v_4$} edge from parent[dashed]}
				}
			}
			child{node {7}}
                child [xshift=-10mm]{node [val] {7\\$\color{brown}v_6$} edge from parent[dashed]}
		}
		child [xshift=-3mm] {node {6} edge from parent[draw=none]
			child [xshift=0mm]{node {7}}
		}
		child [xshift=-20mm]{node {7} edge from parent[draw=none]
            child {node [val] {5\\$\color{brown}v_6$} edge from parent[dashed]}
        }
		child [xshift=-37mm]{node {8} edge from parent[draw=none]};
	\end{tikzpicture}
 \caption{\cmpdformat{} for rank-4 non-uniform hypergraph in \Cref{tab:example}. The special leaves are represented by rectangular nodes with the dropped vertex over the value of the corresponding entry in \blowupT. Note that these special leaves are depicted on different layers of the forest corresponding to the edge size, for instance, the edge $\{1,4,6\}$ corresponds to 3 special leaves all at level 3 (those with value $\color{violet}v_5$, while the edge \{5,7\} correspond to 2 special leaves all at level 3 (those with value $\color{brown}v_6$). }
\label{fig:ccss}
\end{figure*}

\subsection{Space complexity}
The total number of nodes in the \cmpdformat{} depends not only on the number of edges in $\mathcal{E}(\blowupT)$, but also the relative sizes and intersection structure.  Thus, although a single edge $e$ requires at $\binom{|e|}{\ell}$ vertices at level $\ell$ in the forest and $2^{|e|}-1$ vertices in the forest overall, because of the intersections across edges the size of \cmpdformat{} is typically much smaller than the worst case $\sum_{e \in E} 2^{|e|}-1.$

In the next section, we outline how \cmpdformat{} can be used to compute \sttsvc{} on the \blowup{} tensor.
\section{\sttsvc{} computation} \label{sec:s3ttvc}
This work adopts the generating function approach outlined in the AAY algorithm \cite{Aksoy2023ScalableHypergraphs} for computing \sttsvc{} in parallel using the \cmpdformat{} data structure. We present two algorithms -- a baseline approach given in \Cref{alg:naive} which directly parallelizes the AAY approach given in \Cref{alg:without-css} and uses the \cmpdformat{} to store the hypergraph and an optimized version in \Cref{alg:optimized} which rearranges the computations of the AAY approach in order to leverage the \cmpdformat{} to reduce the overall computation by memoization of intermediate results.  

For both of these approaches the focus is on calculating, for every pair $e \in E$ and $v \in e$, the last entry of the convolution in the convolution of lists $E_N(v)$ and $\{ \overline{E}_N(\Vc{b}_u)\}_{u \in e\backslash v}$. For the baseline algorithms, we consider two different methods of computing this convolution; an in-place shift-and-multiply approach and a more efficient approach (but with a larger memory footprint) based on the Fast Fourier Transform (FFT)~\cite{JamesWCooley1965AnSeries}.  In \Cref{alg:optimized},  we only consider a variant of the shift-and-multiply convolution because of the memoization approach used. 

\subsection{Baseline algorithm}
\begin{algorithm}[h]
\small
    \begin{algorithmic}[1]
        \Statex \textbf{Input}: Non-uniform hypergraph stored in \cmpdformat{}, $\V{b}$
        \Statex \textbf{Output}: \sttsvc{} output, $\V{s} = \T{B}\V{b}^{N - 1}$
        \For{$\ell = N, \ldots ,1$} 
            \ParFor{$(e,v) \in \mathcal{L}_{\ell}$}\Comment{\textrm{\cmpdformat{}}}
                    \State $\texttt{coefs} = E_N(\Vc{b}_{v})$
                    \State $u = v$
                    \For{$\ell' = \ell - 1, \ldots , 1$}
                        \State $u$ = \texttt{parent}($u$) \Comment{\textrm{\cmpdformat{}}}
                        \State $\texttt{coefs} = \overline{E}_N(\V{b}_u) \conv \texttt{coefs}$
                    \EndFor
                    \State AtomicAdd$\br*{\V{s}_{v}, \frac{(N-1)!\ell}{|\beta(e)|} \texttt{coefs}[N-1]}$
            \EndParFor
        \EndFor
    \end{algorithmic}
    \caption{\sttsvc{} using \cmpdformat{} via generating function.}
    \label{alg:naive}
\end{algorithm}

Note that \Cref{alg:naive} iterates over the ``special" leaves in the \cmpdformat{} structure by level and for each of these leaves moves up through the \cmpdformat{} forest to a root of one of the subtrees.  This leaf-to-root traversal would make any memoization approach inefficient and challenging to implement as the repeated calculations occur at different locations in the \cmpdformat{}.  For example, in \Cref{fig:ccss} there are several repeated paths (for example, $4,6$ with a special leaf for vertex 1 appears in the trees rooted at 3 and at 4, or 7 with a special leaf for vertex 5 appears in trees rooted at 4 and at 7), however as these calculations occur at different nodes and different levels in the \cmpdformat{} forest it is challenging to realize the benefits of memoization.  This observation inspires our development of a root-to-leaf traversal of the \cmpdformat{} structure, detailed in the next subsection.

\subsection{Optimization: Convolution Memoization}
\begin{algorithm}[h]
    \begin{algorithmic}[1]
        \Statex \textbf{Input}: Non-uniform hypergraph stored in \cmpdformat{}, $\V{b}$
        \Statex \textbf{Output}: \sttsvc{} output, $\V{s} = \blowupT\V{b}^{N- 1}$
        \State For each processor allocate sub-coefficient memoization workspace $W$ of size $\Rn{(N-1) \times (N-1)}$.
        \ParFor{$v  = 1, 2, \ldots n$}
            \For{$j = 1, 2, \ldots N - 1$}
                \For{$i = 1, 2, \ldots, j$}
                    \State $W_{ij} \leftarrow \frac{1}{(j-i+1)!} \V{b}_{v}^{j-i+1}$
                \EndFor
            \EndFor
            \State DFS($v$, W, 1)
        \EndParFor
        \Statex
        \Function{DFS}{$v$, W, $\ell$}
            \For{$u \in \mathcal{S}(v)$}
                \State $\texttt{coefs} = \sqbr*{1, \V{b}_u, \frac{\V{b}_u^2}{2!}, \ldots, \frac{\V{b}_u^{N-1-\ell}}{(N-1-\ell)!}}$
                \State AtomicAdd$\br*{\V{s}_u, \frac{(N-1)!\ell}{\beta(\ell)} \texttt{coefs}^TW_{N-1-\ell}}$
            \EndFor
            \For{$u \in \text{children of }v$}
                \For{$p = 1, 2, \ldots N - 2$}
                    \For{$q = 1, 2, \ldots, p$}
                        \State $Z_{pq} \leftarrow \sum_{c = 0}^{q-p} \frac{1}{c!} \V{b}_u^{c+1} W_{p, q-c}$
                    \EndFor
                \EndFor
                \State DFS($u$, $Z$, $\ell+1$)
            \EndFor
        \EndFunction
    \end{algorithmic}
    \caption{Memoized \sttsvc{} using \cmpdformat{} and generating function.}
    \label{alg:optimized}
\end{algorithm}
In this approach, we optimize the traversal of the \cmpdformat{} forest by using a depth-first search to traverse each tree independently with a separate memoization space, $W$.  After the algorithm has processed a node $v$ on level $\ell$ which has path to the root $v_1, v_2, \ldots, v_{\ell} = v$, the $k^{\textrm{th}}$ column of $W$, denoted $W_k$, stores the portion of the convolution of $\{ \overline{E}_N(\V{b}_{v_i})\}_{i=1}^k$ with degree at most $\ell+k-1$.   Furthermore, if $\mathcal{S}(v)$ is non-empty for a vertex $v$ at level $\ell$, for any $u \in \mathcal{S}(v)$ the contribution of to $\V{s}_u$ can be found by taking the dot product of column $W_{N-1-\ell}$ with the terms from $E_N(\V{b}_u)$ of degree at most $N-1-\ell$.  The pseudocode for this approach is given in \Cref{alg:optimized}.

To illustrate the memoization approach, consider the traversal of the subtree rooted at 1 in Figure 2. This tree will have 3 workspaces associated with it $W_1$, $W_2$, $W_3$ associated with each level of the tree. Workspace 1 will always contain the information necessary to construct the generating function $\bar{E}_N(1)$, while $W_2$ and $W_3$ will contain the information necessary to construct the generating function for the convolutions of $\bar{E}_N$ for $v_1, v_2$ and $v_1, v_2, v_3$, where $v_1, v_2, v_3$ is the path to the current node in the depth-first traversal of the tree. After the unique child of the path (1,2,8) has been computed, the next node in the depth-first traversal is vertex 3 as a child of the root (vertex 1). The updated $W_2$ can be computed directly from the information in $W_1$ while updated $W_3$ is delayed. Now, when traversing the children of vertex 2 (namely 4, 6, and 8) the appropriate $W_3$ can be computed directly from $W_2$ without recomputing the convolutions $\bar{E}_N(1) * \bar{E}_N(2)$. This convolution has been effectively memoized for future computations in $W_2$.

We note that for readability of \Cref{alg:optimized} we have suppressed the use of several easy optimizations.  For instance, if all the edges associated with special leaves in a tree have size at least $k$, then the number of rows $W$ can be reduced to $N - k +1$ as the higher order terms in $\overline{E}_N$ are irrelevant to the final output. Similarly, if the maximum size of an edge associated with a tree is $m$, then the number of columns of $W$ can be reduced to $m-1$ as the longest path to be tracked has size $m$. 
The memoization workspace is allocated per processor, which stores the $W$ matrix of convolution operations. Moreover, note that $W$ is a square upper-triangular Toeplitz matrix, which brings down memory costs to $\mathcal{O}(N)^2$, since each vertex needs to update only a vector of length $N$.
Finally, the \cmpdformat{} forest can be trimmed to eliminate trees, or subtrees, which have no special vertices attached.

\subsection{Computational complexity}
We compare the computational complexity of \Cref{alg:naive} and \Cref{alg:optimized}. For \Cref{alg:naive}, every special leaf corresponding to $(e,v)$ requires one traversal from the special leaf to the root, and thus total number of convolutions computed is $\sum_{e \in E} |e|^2 - |e|$.  However, in \Cref{alg:optimized} every edge of the \cmpdformat{} forest corresponds to a convolution which is computed exactly once.  In particular, the number of convolutions is one less than the number of nodes in the \cmpdformat{} forest.  While the precise speedup resulting from avoiding the extra convolutions is heavily dependent on the structure of the hypergraph, as we will show in \Cref{sec:experiments}, in most real-world cases this results in an order of magnitude saving in runtime.

\begin{figure}
    \centering
    \includegraphics[width=82pt]{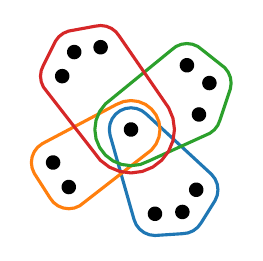}\qquad
    \includegraphics{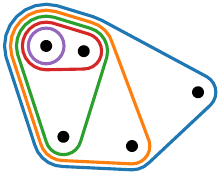}
    \caption{Hypergraph sub-structures in which the memoization of \Cref{alg:optimized} significantly improves the performance over that of \Cref{alg:naive}. On the left is the sunflower hypergraph where many edges share a common intersection, and on the right is the nested hypergraph where edges satisfy a containment relationship.}
    \label{fig:sunflower}
\end{figure}

In spite of the difficulty in determining  the exact number of convolutions saved by using the \cmpdformat{} structure, we can still identify key substructures that lead to significant performance benefits; namely nested families of edges and ``sunflowers"  (collections of edges with a shared intersection), see \Cref{fig:sunflower}.  For example, with a nested series of edges $e_1 \subset e_2 \subset \cdots \subset e_k$, \Cref{alg:naive} will require the computation of $\sum_{i=1}^{k} |e_i|^2 - |e_i|$ convolutions without memoization.  In contrast, the optimal \cmpdformat{} tree will yield require only $\frac{|e_k|^2 - |e_k| -2}{2} + \sum_{i=1}^k |e_i|$ convolutions.  Similarly, if we consider a sunflower consisting of $m$ edges of size $k$ with a common intersection of size $t$, we can see that the non-memoized computation will require $m(k^2-k)$ convolution calculations while the memoized version requires
\[ \frac{(k-t-1)^2 + (k-t-1)}{2} + (t-2)t + (k-t)m + km\]
convolutions.  In fact, even for a single edge of size $k$, the non-memoized computation requires at $k^2 - k$ convolutions while the memoized version will require $\frac{k^2-k}{2}$.  Thus the memoization decreases by a factor of \emph{at least} 2 times the number of convolution operations necessary to evaluate \sttsvc{}.

\section{Experiments}\label{sec:experiments}

We compare the runtime performance and thread scalability of our shared-memory parallel \sttsvc{} algorithm \algo{} against our two baseline approaches --- \naive{} and \fft{} --- for a collection of real-world and synthetic datasets. 

\subsection{Platform and experimental configurations}
Our experiments were conducted on a shared-memory machine with two 64-core AMD Epyc 7713 CPUs at 2.0 GHz and 512GB DDR4 DRAM. 
This work is implemented using C++ and multi-threading parallelized using OpenMP; all numerical operations are performed using double-precision floating point arithmetic and 64-bit unsigned integers. It is compiled using GCC 10.3.0 and Netlib LAPACK 3.8.0 \cite{Anderson1999LAPACKGuide} for linear algebra routines. The polynomial multiplication optimization is performed via full one-dimensional discrete convolution using Fast Fourier Transform (FFT) implementation from FFTW library \cite{Frigo2005TheFFTW3}.

\subsection{Datasets}
We tested with eight real-world datasets and four synthetic datasets for more analysis, shown in \Cref{tab:datasets}. 

\begin{table*}[]
    \centering
    \caption{Summary of datasets. 
    }
    \label{tab:datasets}
    \begin{tabular}{c | c | c | c | c | c | c | p{6cm}}
        Category & Symbol & Dataset & $|V|$ & $|E|$ & $N_k$ & Ref. & Brief overview of dataset\\
        \hline
        \multirow{8}{*}{Real-world} & R1 & MAG-10 & 80198 & 51889 & 25 & \cite{Amburg2020ClusteringLabels,Sinha2015AnApplications} & subset of coauthorship MicrosoftAcademic Graph\\
        & R2 & DAWN & 2109 & 87104 & 22 & \cite{Amburg2020ClusteringLabels} & drugs in patients prior to ER visit\\
        & R3 & cooking & 6714 & 39774 & 65 & \cite{Amburg2020ClusteringLabels} & recipes formed by combining different ingredients\\
        & R4 & walmart-trips & 88860 & 69906 & 25 & \cite{Amburg2020ClusteringLabels} & sets of co-purchased products at Walmart\\
        & R5 & trivago-clicks & 172738 & 233202 & 86 & \cite{Chodrow2021HypergraphModularity} & sets of hotel accommodations clicked out by user\\
        & R6 & amazon-reviews & 2193601 & 3685588 & 26 & \cite{Ni2019JustifyingAspects} & (filtered) sets of Amazon product reviews\\
        & R7 & mathoverflow& 5176 & 39793 & 100 & \cite{Veldt2020MinimizingHypergraphs} & (filtered) sets of answers by Math Overflow users\\
        & R8 & stackoverflow & 7816553 & 1047818 & 76 & \cite{Veldt2020MinimizingHypergraphs} & (filtered) sets of answers by Stack Overflow users\\
        \hline
        \parbox[c]{.4in}{\centering Synthetic} & \parbox[c]{.1in}{\centering S1 \\ S2 \\ S3 \\ S4} & \parbox[c]{.6in}{\centering Synthetic 1 \\ Synthetic 2 \\ Synthetic 3 \\ Synthetic 4} &
        \parbox[c]{.4in}{\centering 2500 \\ 2500 \\ 2500 \\ 2500} & \parbox[c]{.4in}{\centering 25000 \\ 25000 \\ 25000 \\ 25000} & \parbox[c]{.2in}{\centering 30 \\ 40 \\ 50 \\ 60} & \parbox[c]{.2in}{\centering --} & \parbox[c]{2.25in}{\vspace{2pt}Randomly generated non-uniform hypergraphs consisting of 5-, 10-, 15-, ... $N_k-$uniform hypergraphs, with each uniform hypergraph containing $|E|/k$ hyperedges uniformly chosen from all $\binom{|V|}{N_l}$ hyperedges, $1 \leq l \leq k$} \vspace{2pt} \\       
        \hline
    \end{tabular}
\end{table*} 

\emph{Real-world data.} The real-world data come from diverse applications with different amount of nodes $|V|$ (ranging from $10^3$ to $10^6$), hyperedges $|E|$ (ranging from $10^4$ to $10^6$), and component adjacency tensors $N_k$ (ranging from 22--100).
A brief introductory description and reference of the real-world datasets are given in the last two columns. Note that datasets which contain ``filtered'' in their description are filtered versions of the originals, using the less than or equal to filtering from Landry et al.~\cite{Landry2023FilteringDatasets} to remove all hyperedges of size larger than $N_k$.

\emph{Synthetic data.} The synthetic datasets are random hypergraphs on the same set of nodes with the same total amount of hyperedges. But the hyperedges are of varying sizes, uniformly chosen at random over the vertex set. For each of S1, S2, S3 and S4, the component symmetric adjacency tensor orders were taken to be multiples of five until the hypergraph rank i.e. maximum component adjacency tensor order, 
with each component tensor containing approximately the same number of hyperedges.

\subsection{Overall performance}

We evaluate the performance of both implementations of \Cref{alg:naive}, i.e., \naive{} that directly computes the polynomial multiplication and \fft{} that substitutes the polynomial multiplication kernel with a FFT convolution approach, and \Cref{alg:optimized} (\algo{}), which memoizes the coefficients to compute \sttsvc{} using a single traversal of the \cmpdformat{} representation on both categories of datasets. \par

\begin{table}[]
    \centering
    \scriptsize
    \caption{Single-core speedups of \naive{}, \fft{}, and \algo{} versus the Python implemented Algorithm~\ref{alg:without-css} in~\cite{Aksoy2023ScalableHypergraphs}.}
    \begin{tabular}{c|c|c|c|c|c|c|c|c|c}
    speedup  & R1 & R2 & R3 & R4 & R5 & R6 & R7 & R8  \\
                \hline
        \naive{} & 0.7 & 1.0 & 0.2 & 1.4 & 0.1 & - & 0.1 & - \\
        \fft{}   & 2.4 & 3.5 & 2.8 & 4.2 & 0.4 & 4.7 & 2.1 & 2.6 \\ 
        \algo{}  & 2.4 & 5.2 & 3.0 & 17.0 & 1.3 & 18.6 & 4.5 & 4.6 \\
    \end{tabular}
    \label{tab:speedup_vs_python}
\end{table}

\textbf{Comparisons to SOTA methods.}
In Table~\ref{tab:speedup_vs_python}, we summarize the speedups of \naive{}, \fft{}, and \algo{} on a single core as compared to the state-of-the-art single-core Python implementation of Algorithm~\ref{alg:without-css} \cite{Aksoy2023ScalableHypergraphs}. 
\footnote{In personal conversations with the authors of~\cite{Aksoy2023ScalableHypergraphs}, it was noted that a highly-optimized low-level FFT implementation with a Python front-end was used to obtain their results.}
Our \algo{} outperforms considerably, providing speedups of 1.3 -- 18.6 times on the eight real-world datasets, while \fft{} provides considerable improvements as well. \naive{} largely underperforms, due to using direct convolution instead of FFT, with the exception being \texttt{walmart-trips}, where it obtains a speedup of 1.4. 
Although the timing results in Aksoy et al.~\cite{Aksoy2023ScalableHypergraphs} were obtained using a different language and system, the speedup observed is almost two orders of magnitude for the largest dataset, \texttt{amazon-reviews}, which suggests that
 the speedup is due to improvements in the algorithm as opposed to system and implementation differences. 
Further, as the dominate subroutine (the convolution operation) is implemented using optimized and compiled code (via a Python-wrapper in \cite{Aksoy2023ScalableHypergraphs}) one would expect that much of the performance differences resulting from language choice are mitigated. Furthermore, the strong scaling results seen in~\Cref{fig:s3ttvc-scalability} suggest even greater speedups as we increase the number of cores.

\begin{figure}
    \centering
    \includegraphics[scale=0.5]{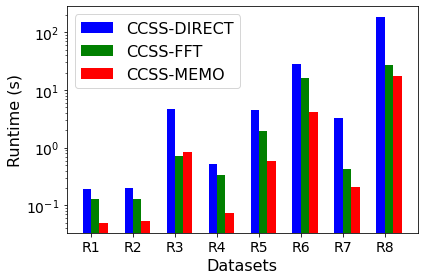}
    \caption{Overall runtime performance of \algo{}, \naive{} and \fft{} for real-world datasets in \Cref{tab:datasets}}
    \label{fig:s3ttvc-runtime}
\end{figure}

\textbf{Real-world datasets.} \Cref{fig:s3ttvc-runtime} presents the runtime performance on 128 cores of these three approaches for the real-world datasets. \algo{} performs the best on almost all the eight datasets,  achieves $1.97-53.98\times$ 
speedup over \naive{} and $0.65-12.45\times$ 
speedup over \fft{}. The variation in speedups across different datasets for \algo{} is inherent to the hypergraph structure, and is due to the degree of overlap in the hyperedges present in the hypergraph. For the R3 dataset, we can see from \Cref{fig:ccss-construction} that \cmpdformat{} achieves very low compression compared to the coordinate format. Thus, there is no significant advantage in using the memoization approach to compute \sttsvc{}, which is why we see that \fft{} slightly outperforms \algo{} for higher thread configurations for this dataset. 

\begin{figure}
    \centering
    \includegraphics[scale=0.5]{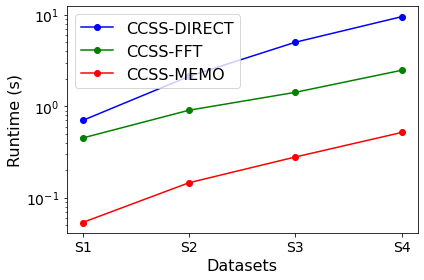}
    \caption{Comparison of runtime of \algo{} with \naive{} and \fft{} for synthetic datasets in \Cref{tab:datasets} to highlight effect of rank of non-uniform hypergraph on \sttsvc{} computation. }
    \label{fig:s3ttvc-memoization}
\end{figure}

\textbf{Synthetic datasets.} Since the synthetic datasets maintain the same number of IOU non-zeros across component uniform hypergraphs, the number of leaf nodes per level of the \cmpdformat{} that contribute the \sttsvc{} computation is the same. This allows us to inspect the effect of the rank of the non-uniform hypergraph on the performance of all three algorithms. We see from \Cref{fig:s3ttvc-memoization} that sub-coefficient memoization has a significant impact on performance for hypergraphs of increasing ranks. This is to be expected since with increasing tensor order, the reduction in the number of traversals of the \cmpdformat{}, as well as sharing of sub-coefficients between overlapping hyperedges for a uniform hyperedge distribution, would result in improved performance compared to \naive{} and \fft{}.
\par

\begin{figure*}[h]
    \centering
    \includegraphics[scale=0.5]{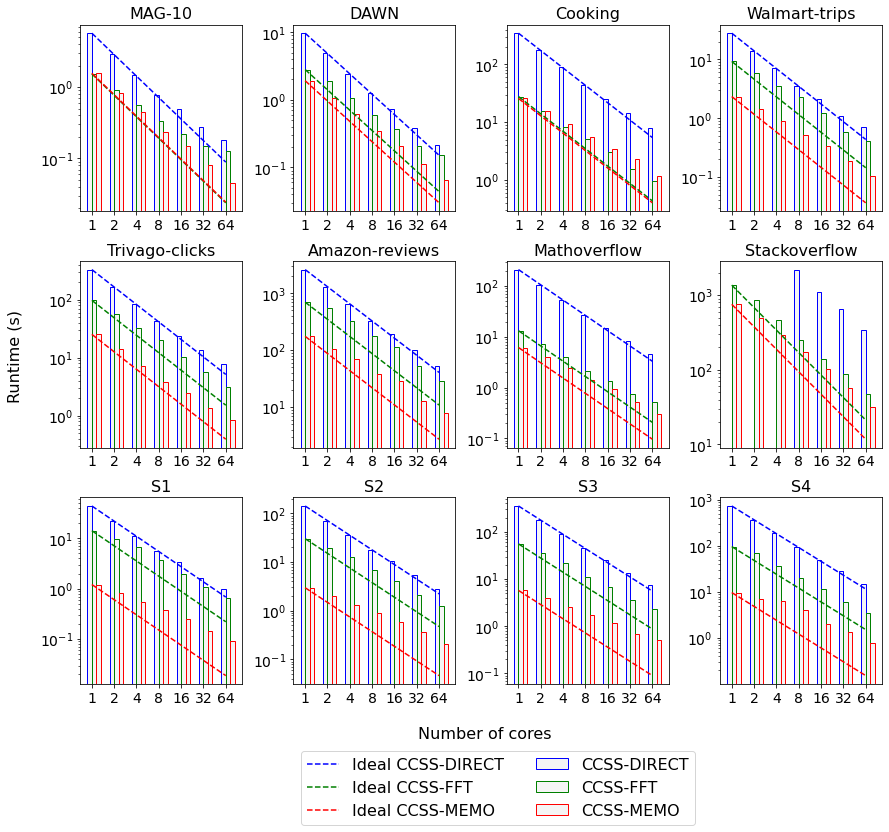}
    \caption{Thread scalability of \naive{}, \fft{}, and \algo{} for synthetic and real-world datasets in \Cref{tab:datasets}. We set the time limit of our jobs at 50 minutes for all thread configurations. The \texttt{stackoverflow} dataset does not complete within the allotted time for 1, 2 and 4 thread configurations for \naive{}.}
    \label{fig:s3ttvc-scalability}
\end{figure*}

\subsection{Thread scalability}
\Cref{fig:s3ttvc-scalability} compares the thread scalability of \naive{}, \fft{}, and \algo{} for synthetic and real-world datasets in \Cref{tab:datasets}. Across all datasets, we see that the \algo{} approach is faster than both the baseline approaches. The dashed lines indicate the ideal speedup lines for each of the three algorithms.
\naive{} and \fft{} do not use memoization for the intermediate $\bar{E}_N(v)$ computations.
The \naive{} algorithm shows the best scalability of the three approaches, while both \fft{} and \algo{} show decreasing scalability with increasing number of threads. For both of these approaches, as the work done per thread reduces (in terms of optimized FFT subroutines in FFTW for \fft{} and in terms of coefficient $W$ memoization for \algo{}), the overhead in the atomic operation becomes more significant, especially for the hypergraphs with smaller number of nodes, and this manifests as suboptimal scaling. 
Moreover, the scalability of \algo{} is also affected by the load imbalance between threads due to the varying number of IOU non-zeros across trees within the CCSS structure.

\begin{figure}
    \centering
    \includegraphics[scale=0.45]{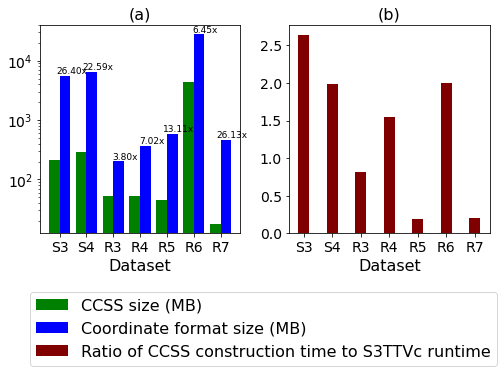}
    \caption{Summary statistics for \cmpdformat{} construction for representative datasets in \Cref{tab:datasets}. We see that \cmpdformat{} achieves more compact storage than the coordinate format i.e. storing the adjacency tensor non-zeros as lists of index tuples along with their non-zero value. }
    \label{fig:ccss-construction}
\end{figure}

\subsection{\cmpdformat{} Construction}
\textbf{\cmpdformat{} construction.} While the \format{} compressed only in terms of overlapping IOU non-zeros within a symmetric adjacency tensor, the \cmpdformat{} adds another layer of compactness in terms of shared indices between IOU non-zeros across multiple symmetric tensors. Moreover, for computing \sttsvc{} on the \blowup{} tensor using the \cmpdformat{}, we can further prune paths if the leaf node of the path does not own any IOU non-zeros. \Cref{fig:ccss-construction}(a) shows the size of the constructed \cmpdformat{} for representative synthetic and real-world datasets, while \Cref{fig:ccss-construction}(b) examines the ratio of the amount of time spent in the construction of \cmpdformat{} to the runtime of \algo{} for \sttsvc{} computation.
\par

\subsection{$H$-eigenvector centrality computation speedups}
 \algo{} provides a practical framework to compute centrality in large hypergraphs. We compute tensor H-eigenvector centrality using \Cref{alg:nqz} on the largest real-world dataset in \Cref{tab:datasets} --- the \texttt{amazon-reviews} --- on 128 cores. \algo{} obtains speedups of $6.49\times$ and $3.53\times$ over \naive{} and \fft{} respectively, which shows the applicability of this work in the analysis of real-world hypergraphs.
\section{Conclusion} \label{sec:conclusion}
This work introduces the \cmpdformat{}, an extension of the \format{} for uniform hypergraphs, to compactly non-uniform hypergraphs. A novel memoization-based algorithm \algo{} adapted from the generating function approach is proposed to compute \sttsvc{} on the \blowup{} adjacency tensor of non-uniform hypergraphs. We demonstrate the performance of our shared-memory parallel \algo{} by comparing it to two naive baseline algorithms using the \cmpdformat{} - \naive{} and \fft{} - for multiple synthetic and real-world datasets. In the future, we plan to explore distributed-memory construction of \cmpdformat{} and distributed-memory parallel \sttsvc{} computation using \algo{}. Furthermore, we plan to utilize \algo{} as the computational kernel in extending multilinear PageRank to nonuniform hypergraphs, where we anticipate its advanced data structures and parallelization will allow analysis of datasets previously considered prohibitively large for tensor analysis. Our fast algorithm would also help facilitate the multilinear hypergraph clustering ~\cite{Aksoy2023ScalableHypergraphs}. Moreover, we believe it opens the door for development of tensor-based methods for semi-supervised and supervised hypergraph learning tasks such as node classification and link prediction.

\bibliographystyle{IEEEtran}
\bibliography{references}

\end{document}